\documentclass[a4paper,12pt]{article}
\usepackage[all]{xy}
\usepackage[T1]{fontenc}
\usepackage[dvips]{graphicx}
\usepackage{amsmath,amssymb,amsfonts,amsthm,latexsym,mathrsfs,textcomp,verbatim}
\xyoption{web}

\title{One topos, many sites}
\author{{Olivia Caramello} \vspace{3 mm}\\ {\small DPMMS, University of Cambridge,}\\{\small Wilberforce Road, Cambridge CB3 0WB, UK}\\{\small O.Caramello@dpmms.cam.ac.uk}}
\date{July 14, 2009}
\begin{document}
\bgroup           
\let\footnoterule\relax  
\maketitle  
\begin{abstract}
We give characterizations, for various fragments of geometric logic, of the class of theories classified by a locally connected (resp. connected and locally connected, atomic, compact, presheaf) topos, and exploit the existence of multiple sites of definition for a given topos to establish some properties of quotients of theories of presheaf type.  
\end{abstract} 
\egroup 
\vspace{5 mm}


\def\Monthnameof#1{\ifcase#1\or
   January\or February\or March\or April\or May\or June\or
   July\or August\or September\or October\or November\or December\fi}
\def\today{\number\day~\Monthnameof\month~\number\year}

%
%
%
\def\pushright#1{{
   \parfillskip=0pt            
   \widowpenalty=10000         
   \displaywidowpenalty=10000  
   \finalhyphendemerits=0      
  %
   \leavevmode                 
   \unskip                     
   \nobreak                    
   \hfil                       
   \penalty50                  
   \hskip.2em                  
   \null                       
   \hfill                      
   {#1}                        
  %
   \par}}                      

\def\qed{\pushright{$\square$}\penalty-700 \smallskip}

\newtheorem{theorem}{Theorem}[section]

\newtheorem{proposition}[theorem]{Proposition}

\newtheorem{scholium}[theorem]{Scholium}

\newtheorem{lemma}[theorem]{Lemma}

\newtheorem{corollary}[theorem]{Corollary}

\newtheorem{conjecture}[theorem]{Conjecture}

\newenvironment{proofs}%
 {\begin{trivlist}\item[]{\bf Proof }}%
 {\qed\end{trivlist}}

  \newtheorem{rmk}[theorem]{Remark}
\newenvironment{remark}{\begin{rmk}\em}{\end{rmk}}

  \newtheorem{rmks}[theorem]{Remarks}
\newenvironment{remarks}{\begin{rmks}\em}{\end{rmks}}

  \newtheorem{defn}[theorem]{Definition}
\newenvironment{definition}{\begin{defn}\em}{\end{defn}}

  \newtheorem{eg}[theorem]{Example}
\newenvironment{example}{\begin{eg}\em}{\end{eg}}

  \newtheorem{egs}[theorem]{Examples}
\newenvironment{examples}{\begin{egs}\em}{\end{egs}}


\mathcode`\<="4268  
\mathcode`\>="5269  
\mathcode`\.="313A  
\mathchardef\semicolon="603B 
\mathchardef\gt="313E
\mathchardef\lt="313C

\newcommand{\app}
 {{\sf app}}

\newcommand{\Ass}
 {{\bf Ass}}

\newcommand{\ASS}
 {{\mathbb A}{\sf ss}}

\newcommand{\Bb}
{\mathbb}

\newcommand{\biimp}
 {\!\Leftrightarrow\!}

\newcommand{\bim}
 {\rightarrowtail\kern-1em\twoheadrightarrow}

\newcommand{\bjg}
 {\mathrel{{\dashv}\,{\vdash}}}

\newcommand{\bstp}[3]
 {\mbox{$#1\! : #2 \bim #3$}}

\newcommand{\cat}
 {\!\mbox{\t{\ }}}

\newcommand{\cinf}
 {C^{\infty}}

\newcommand{\cinfrg}
 {\cinf\hy{\bf Rng}}

\newcommand{\cocomma}[2]
 {\mbox{$(#1\!\uparrow\!#2)$}}

\newcommand{\cod}
 {{\rm cod}}

\newcommand{\comma}[2]
 {\mbox{$(#1\!\downarrow\!#2)$}}

\newcommand{\comp}
 {\circ}

\newcommand{\cons}
 {{\sf cons}}

\newcommand{\Cont}
 {{\bf Cont}}

\newcommand{\ContE}
 {{\bf Cont}_{\cal E}}

\newcommand{\ContS}
 {{\bf Cont}_{\cal S}}

\newcommand{\cover}
 {-\!\!\triangleright\,}

\newcommand{\cstp}[3]
 {\mbox{$#1\! : #2 \cover #3$}}

\newcommand{\Dec}
 {{\rm Dec}}

\newcommand{\DEC}
 {{\mathbb D}{\sf ec}}

\newcommand{\den}[1]
 {[\![#1]\!]}

\newcommand{\Desc}
 {{\bf Desc}}

\newcommand{\dom}
 {{\rm dom}}

\newcommand{\Eff}
 {{\bf Eff}}

\newcommand{\EFF}
 {{\mathbb E}{\sf ff}}

\newcommand{\empstg}
 {[\,]}

\newcommand{\epi}
 {\twoheadrightarrow}

\newcommand{\estp}[3]
 {\mbox{$#1 \! : #2 \epi #3$}}

\newcommand{\ev}
 {{\rm ev}}

\newcommand{\Ext}
 {{\rm Ext}}

\newcommand{\fr}
 {\sf}

\newcommand{\fst}
 {{\sf fst}}

\newcommand{\fun}[2]
 {\mbox{$[#1\!\to\!#2]$}}

\newcommand{\funs}[2]
 {[#1\!\to\!#2]}

\newcommand{\Gl}
 {{\bf Gl}}

\newcommand{\hash}
 {\,\#\,}

\newcommand{\hy}
 {\mbox{-}}

\newcommand{\im}
 {{\rm im}}

\newcommand{\imp}
 {\!\Rightarrow\!}

\newcommand{\Ind}[1]
 {{\rm Ind}\hy #1}

\newcommand{\iten}[1]
{\item[{\rm (#1)}]}

\newcommand{\iter}
 {{\sf iter}}

\newcommand{\Kalg}
 {K\hy{\bf Alg}}

\newcommand{\llim}
 {{\mbox{$\lower.95ex\hbox{{\rm lim}}$}\atop{\scriptstyle
{\leftarrow}}}{}}

\newcommand{\llimd}
 {\lower0.37ex\hbox{$\pile{\lim \\ {\scriptstyle
\leftarrow}}$}{}}

\newcommand{\Mf}
 {{\bf Mf}}

\newcommand{\Mod}
 {{\bf Mod}}

\newcommand{\MOD}
{{\mathbb M}{\sf od}}

\newcommand{\mono}
 {\rightarrowtail}

\newcommand{\mor}
 {{\rm mor}}

\newcommand{\mstp}[3]
 {\mbox{$#1\! : #2 \mono #3$}}

\newcommand{\Mu}
 {{\rm M}}

\newcommand{\name}[1]
 {\mbox{$\ulcorner #1 \urcorner$}}

\newcommand{\names}[1]
 {\mbox{$\ulcorner$} #1 \mbox{$\urcorner$}}

\newcommand{\nml}
 {\triangleleft}

\newcommand{\ob}
 {{\rm ob}}

\newcommand{\op}
 {^{\rm op}}

\newcommand{\pepi}
 {\rightharpoondown\kern-0.9em\rightharpoondown}

\newcommand{\pmap}
 {\rightharpoondown}

\newcommand{\Pos}
 {{\bf Pos}}

\newcommand{\prarr}
 {\rightrightarrows}

\newcommand{\princfil}[1]
 {\mbox{$\uparrow\!(#1)$}}

\newcommand{\princid}[1]
 {\mbox{$\downarrow\!(#1)$}}

\newcommand{\prstp}[3]
 {\mbox{$#1\! : #2 \prarr #3$}}

\newcommand{\pstp}[3]
 {\mbox{$#1\! : #2 \pmap #3$}}

\newcommand{\relarr}
 {\looparrowright}

\newcommand{\rlim}
 {{\mbox{$\lower.95ex\hbox{{\rm lim}}$}\atop{\scriptstyle
{\rightarrow}}}{}}

\newcommand{\rlimd}
 {\lower0.37ex\hbox{$\pile{\lim \\ {\scriptstyle
\rightarrow}}$}{}}

\newcommand{\rstp}[3]
 {\mbox{$#1\! : #2 \relarr #3$}}

\newcommand{\scn}
 {{\bf scn}}

\newcommand{\scnS}
 {{\bf scn}_{\cal S}}

\newcommand{\semid}
 {\rtimes}

\newcommand{\Sep}
 {{\bf Sep}}

\newcommand{\sep}
 {{\bf sep}}

\newcommand{\Set}
 {{\bf Set }}

\newcommand{\Sh}
 {{\bf Sh}}

\newcommand{\ShE}
 {{\bf Sh}_{\cal E}}

\newcommand{\ShS}
 {{\bf Sh}_{\cal S}}

\newcommand{\sh}
 {{\bf sh}}

\newcommand{\Simp}
 {{\bf \Delta}}

\newcommand{\snd}
 {{\sf snd}}

\newcommand{\stg}[1]
 {\vec{#1}}

\newcommand{\stp}[3]
 {\mbox{$#1\! : #2 \to #3$}}

\newcommand{\Sub}
 {{\rm Sub}}

\newcommand{\SUB}
 {{\mathbb S}{\sf ub}}

\newcommand{\tbel}
 {\prec\!\prec}

\newcommand{\tic}[2]
 {\mbox{$#1\!.\!#2$}}

\newcommand{\tp}
 {\!:}

\newcommand{\tps}
 {:}

\newcommand{\tsub}
 {\pile{\lower0.5ex\hbox{.} \\ -}}

\newcommand{\wavy}
 {\leadsto}

\newcommand{\wavydown}
 {\,{\mbox{\raise.2ex\hbox{\hbox{$\wr$}
\kern-.73em{\lower.5ex\hbox{$\scriptstyle{\vee}$}}}}}\,}

\newcommand{\wbel}
 {\lt\!\lt}

\newcommand{\wstp}[3]
 {\mbox{$#1\!: #2 \wavy #3$}}

\newcommand{\fu}[2]
{[#1,#2]}

\newcommand{\st}[2]
 {\mbox{$#1 \to #2$}}

\tableofcontents

\newpage
\section{Introduction}

Given the fact that Grothendieck toposes are `the same thing as' Morita-equivalence classes of geometric theories, it is naturally of interest to investigate how classical topos-theoretic properties of toposes translate into logical properties of the theories they classify. 

Characterizations of the class of geometric theories classified by a Boolean (resp. De Morgan) topos have been provided in \cite{OC3}. In the third section of this paper, we provide syntactic characterizations, for various fragments of geometric logic, of the class of theories classified by a locally connected (resp. connected and locally connected, atomic, compact, presheaf) topos. Also, we establish criteria for a geometric theory over a given signature to be cartesian (resp. regular, coherent). 

In the last section, given a quotient ${\mathbb T}'$ of a theory of presheaf type $\mathbb T$ correponding to a Grothendieck topology $J$ on the opposite of the category of finitely presentable $\mathbb T$-models as in \cite{OC6}, we discuss how `geometrical' properties of $J$ translate into syntactic properties of ${\mathbb T}'$. In this context, we also show that, given a theory of presheaf type $\mathbb T$, the category of finitely presentable models of $\mathbb T$ is equivalent to a full subcategory of the syntactic category of $\mathbb T$. 

Our technique to establish these latter results is to transfer topos-theoret-\\ic invariants (i.e. properties of objects of toposes which are stable under topos-theoretic equivalence) from one site of definition of a given topos to another.   

The terminology used in the paper is borrowed from \cite{El} and \cite{El2}, if not otherwise stated. Our notion of site is that that of a small Grothendieck site; anyway, all the results on sites established in the paper can be (trivially) extended to essentially small sites of the form $({\cal C}, J)$ where $\cal C$ is an essentially small category.  

\section{Some geometric invariants}

Let us start with some general facts about dense subcategories. We recall from \cite{El2} the following definition.
\begin{definition} 
Let $({\cal C}, J)$ be a site. We say a subcategory $\cal D$ of $\cal C$ is $J$-dense if

(i) every object $c$ of $\cal C$ has a covering sieve $R \in J(c)$ generated by morphisms whose domains are in $\cal D$; and
 
(ii) for any morphism $f:c \to d$ in $\cal C$ with $d\in {\cal D}$, there is a covering sieve $R \in J(c)$ generated by morphisms $g: b \to c$ for which the composite $f\circ g$ is in $\cal D$.
\end{definition} 

Let us denote by $\mathfrak{Groth}({\cal C})$ the Heyting algebra of Grothendieck topologies on a category $\cal C$ (cfr. \cite{OC6}). Given a subcategory ${\cal C}'$ of $\cal C$, we denote by $\mathfrak{Groth}_{{\cal C}'}({\cal C})$ the subset of $\mathfrak{Groth}({\cal C})$ formed by the Grothendieck topologies $J$ on $\cal C$ such that ${\cal C}'$ is $J$-dense.
 
There is an obvious notion of intersection of subcategories; specifically, given a collection $\{{\cal C}_{i} \textrm{ | } i\in I\}$ of subcategories of a category $\cal C$, we can define their intersection ${\cal C}'$ by putting $ob({{\cal C}'})=\mathbin{\mathop{\textrm{\huge $\cap$}}\limits_{i\in I}} ob({\cal C}_{i})$ and $arr({{\cal C}'})=\mathbin{\mathop{\textrm{\huge $\cap$}}\limits_{i\in I}} arr({\cal C}_{i})$ i.e. given an arrow $f:a\to b$ in $\cal C$, $f$ belongs to $arr({{\cal C}'})$ if and only if it belongs to $arr({{\cal C}_{i}})$ for every $i\in I$. It is immediate to see that ${\cal C}'$ is a subcategory of $\cal C$.

We note that, for any Grothendieck topology $J$ on $\cal C$, any finite intersection of subcategories of $\cal C$ which are $J$-dense is again $J$-dense; indeed, this easily follows from the fact that a finite intersection of $J$-covering sieves is again $J$-covering.

The following result provides a couple of useful facts about dense subcategories.

\begin{proposition}\label{Onetoposmany_denseness}
Let $\cal C$ be a category. Then

(i) $\mathfrak{Groth}_{{\cal C}'}({\cal C})$ is closed in $\mathfrak{Groth}({\cal C})$ under arbitrary (non-empty) intersections (i.e. meets in $\mathfrak{Groth}({\cal C})$) and under taking larger topologies, and hence it is an Heyting algebra inheriting the Heyting operations from $\mathfrak{Groth}({\cal C})$ whose maximal element is the maximal Grothendieck topology on $\cal C$ and minimal element is the intersection of all the topologies in $\mathfrak{Groth}_{{\cal C}'}({\cal C})$; 

(ii) Let ${\cal C}'$ and ${\cal C}''$ subcategories of $\cal C$ such that ${\cal C}'$ is a subcategory of ${\cal C}''$. Then ${\cal C}'$ is $J$-dense if and only if ${\cal C}'$ is $J|_{{\cal C}''}$-dense (as a subcategory of ${\cal C}''$) and ${\cal C}''$ is $J$-dense.    
\end{proposition}
     
\begin{proofs}
(i) The first assertion is obvious while the second easily follows from the fact that arbitrary (non-empty) unions of $J$-covering sieves are $J$-covering.

(ii) Suppose that ${\cal C}'$ is $J$-dense. Then, by part (i) of the proposition, ${\cal C}''$ is $J$-dense, and it is immediate to see that ${\cal C}'$ is $J|_{{\cal C}''}$-dense (as a subcategory of ${\cal C}''$). The converse follows from the fact that `composition' of $J$-covering sieves (in the sense of Definition 2.3 \cite{OC6}) is $J$-covering.
\end{proofs}

We note that a small full subcategory of a Grothendieck topos $\cal E$ is $J_{\cal E}$-dense, where $J_{\cal E}$ is the canonical Grothendieck topology on $\cal E$, if and only if it is a separating set for $\cal E$.

Let us also recall the following notions. 

\begin{definition}
Let $\cal E$ be a Grothendieck topos and $A$ an object of $\cal E$. Then

(i) $A$ is said to be an atom if the only subobjects of $A$ in $\cal E$ are the identity on $A$ and the zero subobject, and they are distinct from each other;

(ii) $A$ is said to be indecomposable if does not admit any non-trivial coproduct decompositions;

(iii) $A$ is said to be irreducible if it is $J_{\cal E}$-irreducible, where $J_{\cal E}$ is the canonical topology on $\cal E$; in other words, if any sieve in $\cal E$ containing a small epimorphic family contains the identity on $A$;

(iv) $A$ is said to be compact if every small covering family $\{A_{i} \to A \textrm{ | } i\in I \}$ contains a finite covering subfamily;

(v) $A$ is said to be coherent if it is compact and, whenever we are given a morphism $f:B\to A$ with $B$ compact, the domain of the kernel-pair of $f$ is compact;

(vi) $A$ is said to be supercompact if every small covering family $\{A_{i}\to A \textrm{ | } i\in I \}$ contains a cover;

(vii) $A$ is said to be regular if it is supercompact and, whenever we are given a morphism $f:B\to A$ with $B$ supercompact, the domain of the kernel-pair of $f$ is supercompact.
\end{definition}
 
Recall that an object in a locally connected topos is indecomposable if and only if it is connected (cfr. the discussion after the proof of Lemma C3.3.6 \cite{El2}). 

\begin{rmk}\label{Onetoposmany_Rem}
\emph{It readily follows from the definitions that every coherent (resp. regular) object is compact (resp. supercompact), every atom is an indecomposable object, every irreducible object is supercompact, every supercompact object is indecomposable.}
\end{rmk}

Let us recall from \cite{El2} the following terminology.

A site $({\cal C}, J)$ is said to be locally connected if every $J$-covering sieve is connected i.e. for any $R\in J(c)$, $R$ is connected as a full subcategory of ${\cal C}\slash c$.

A site $({\cal C}, J)$ is said to be atomic if $\cal C$ satisfies the right Ore condition and $J$ is the atomic topology on $\cal C$.

Given a site  $({\cal C}, J)$,

(i) we say an object $c$ of $\cal C$ is $J$-irreducible if the only $J$-covering sieve on $c$ is the maximal sieve $M(c)$;

(ii) we say $J$ is rigid if, for every $c\in {\cal C}$, the family of all morphisms from $J$-irreducible objects to $c$ generates a $J$-covering sieve. 
The site $({\cal C}, J)$ is said to be rigid if $J$ is rigid as a Grothendieck topology on $\cal C$.

We will say that a site $({\cal C}, J)$ is coherent (resp. regular) if $\cal C$ is cartesian and $J$ is a finite-type Grothendieck topology on $\cal C$ (resp. a Grothendieck topology $J$ such that every $J$-covering sieve is generated by a single arrow).

Given a site $({\cal C}, J)$, we denote by $l^{\cal C}_{J}:{\cal C}\to \Sh({\cal C}, J)$ the composite of the Yoneda embedding ${\cal C}\to [{\cal C}^{\textrm{op}}, \Set]$ with the associated sheaf functor $a_{J}:[{\cal C}^{\textrm{op}}, \Set] \to \Sh({\cal C}, J)$. We have the following result.

\begin{proposition}\label{Onetoposmany_identif}
Let $({\cal C}, J)$ be a site. Then

(i) if $({\cal C}, J)$ is locally connected then for each $c\in {\cal C}$, $l^{\cal C}_{J}(c)$ is an indecomposable (equivalently, connected) object in $\Sh({\cal C}, J)$;
 
(ii) if $({\cal C}, J)$ is atomic then for each $c\in {\cal C}$, $l^{\cal C}_{J}(c)$ is an atom in $\Sh({\cal C}, J)$;

(iii) if $({\cal C}, J)$ is rigid then for each $c\in {\cal C}$ such that $c$ is $J$-irreducible, $l^{\cal C}_{J}(c)$ is an indecomposable projective (equivalently, an irreducible object) of $\Sh({\cal C}, J)$;

(iv) if $({\cal C}, J)$ is coherent (resp. regular) then for each $c\in {\cal C}$, $l^{\cal C}_{J}(c)$ is a coherent (resp. regular) object of $\Sh({\cal C}, J)$.
\end{proposition}  

\begin{proofs}
(i) and (ii) were proved in \cite{OC5}. 

To prove (iii), we note that if $({\cal C}, J)$ is rigid then the Comparison Lemma yields $\Sh({\cal C}, J) \simeq [{\cal D}^{\textrm{op}}, \Set]$ where $\cal D$ is the full subcategory of $\cal C$ on the $J$-irreducible objects (cfr. the discussion after Definition C2.2.18 \cite{El2}) and, under this equivalence, for each $c\in {\cal D}$, $l^{\cal C}_{J}(d)$ corresponds to the representable $y(d):{\cal D}^{\op}\to \Set$. Now, it is well-known that all the representables on ${\cal D}$ are indecomposable projective objects in $[{\cal D}^{\textrm{op}}, \Set]$, from which our thesis follows.

Part (iv) was proved in \cite{El2} (cfr. Theorem D3.3.7 and Remark D3.3.10).       
\end{proofs}

We note that if $\cal E$ is equivalent to a presheaf topos then any object of $\cal E$ is irreducible if and only if it is indecomposable and projective. Indeed, by the argument in the proof of Lemma C2.2.20 \cite{El2}, any indecomposable projective object in $\cal E$ is irreducible. Conversely, if $\cal G$ is the full subcategory of $\cal E$ on the irreducible objects then, by Lemma C2.2.20, $\cal G$ is essentially small and $({\cal G}, J_{\cal E}|_{\cal G})$ is a rigid site; moreover, by the Comparison Lemma ${\cal E}\simeq \Sh({\cal G}, J_{\cal E}|_{\cal G})$ and hence Proposition \ref{Onetoposmany_identif}(iii) implies our thesis.

Let us recall from \cite{El2} (C2.3.2(c)) that, for any essentially small site $({\cal C}, J)$, a sieve $R$ on $l^{\cal C}_{J}(c)$ in $\Sh({\cal C}, J)$ is epimorphic iff the sieve $\{f:d\to c  \textrm{ | } l^{\cal C}_{J}(f)\in R \}$ is $J$-covering in $\cal C$. This fact enables us to express properties of objects of the form $l^{\cal C}_{J}(c)$ like compactness, supercompactness or irreducibility in terms of properties of $J$-covering sieves on $c$, as follows (point (i) of the following proposition was proved in \cite{El2} as Lemma D3.3.4).

\begin{proposition}\label{Onetoposmany_express}
Let $({\cal C}, J)$ be a site. Then, with the notation above, we have:

(i) $l^{\cal C}_{J}(c)$ is compact if and only if every $J$-covering sieve on $c$ contains a finite family of arrows which generates a $J$-covering sieve;
 
(ii) $l^{\cal C}_{J}(c)$ is supercompact if and only if every $J$-covering sieve on $c$ contains a single arrow which generates a 
$J$-covering sieve;

(iii) $l^{\cal C}_{J}(c)$ is irreducible if and only if every $J$-covering sieve on $c$ is maximal i.e. $c$ is $J$-irreducible;
\end{proposition}\qed

It turns out that one can rephrase many interesting properties of Grothen-\\dieck toposes in terms of the existence of separating sets for them with particular properties. For example, it is well-known (cfr. Lemma C2.2.20 \cite{El2}) that a Grothendieck topos is equivalent to a presheaf topos if and only if it has a separating set of indecomposable projective objects; moreover, we have the following characterizations. 

\begin{theorem}\label{Onetoposmany_car}
Let $\cal E$ be a Grothendieck topos. Then

(i) $\cal E$ is locally connected if and only if it has a separating set of indecomposable objects;

(ii) $\cal E$ is connected and locally connected if and only if it has a separating set of indecomposable objects containing the terminal object of $\cal E$;

(iii) $\cal E$ is atomic if and only if it has a separating set of atoms;

(iv) $\cal E$ is coherent (resp. regular) if and only if it has a separating set of coherent (resp. regular) objects which is closed in $\cal E$ under finite limits.
\end{theorem}  

\begin{proofs}
(i) Suppose that $\cal E$ is locally connected; then, by Theorem C3.3.10 \cite{El2}, $\cal E$ is of the form $\Sh({\cal C}, J)$ for a locally connected small site $({\cal C}, J)$. Then the objects of the form $l^{\cal C}_{J}(c)$ for $c\in {\cal C}$ are indecomposable (by Proposition \ref{Onetoposmany_identif}(i)) and hence they form a separating set for $\cal E$. Conversely, suppose that $\cal E$ has a separating set of indecomposable objects; then, by arguing as in the proof of Theorem C3.3.10 \cite{El2}, we obtain that the full subcategory ${\cal I}$ of $\cal E$ on the indecomposable objects, equipped with the Grothendieck topology $J_{\cal E}|_{\cal I}$ on $\cal I$ induced by the canonical coverage on $\cal E$ is a locally connected site; but, by the Comparison Lemma, $\cal E$ is equivalent to $\Sh({\cal I}, J_{\cal E}|_{\cal I})$, so that the thesis follows from Theorem C3.3.10 \cite{El2}.

(ii) This follows analogously to part (i), by using the `connected and locally connected' version of Theorem C3.3.10 \cite{El2}.

(iii) This was already proved in \cite{OC5} (Proposition 1.3(i)).

(iv)  One direction follows from Theorem D3.3.7 \cite{El2} and Remarks D3.3.9 and D3.3.10 \cite{El2}. Conversely, suppose that $\cal E$ has a separating set $\cal G$ of coherent (resp. regular) objects which is closed in $\cal E$ under finite limits. Then $\cal G$ is a cartesian category and, by the Comparison Lemma, $\cal E$ is equivalent to $\Sh({\cal G}, J_{\cal E}|_{\cal G})$; now, by Proposition \ref{Onetoposmany_express}(i) (resp. Proposition \ref{Onetoposmany_express}(ii)), the site $({\cal G}, J_{\cal E}|_{\cal G})$ is coherent (resp. regular), and hence, by Theorem D3.3.1 \cite{El2}, $\cal E$ is a coherent (resp. regular) topos, as required.   
\end{proofs}

\begin{rmk}
\emph{Notice that it follows from the theorem and Remark \ref{Onetoposmany_Rem} that any presheaf topos is locally connected, any regular topos is locally connected and any atomic topos is locally connected. Moreover, it is clear that if $\cal E$ is a Boolean topos then for any object $A$ of $\cal E$, $A$ is an atom if and only if it is non-zero and indecomposable, from which it follows that any Boolean locally connected topos is atomic.}
\end{rmk}  

\begin{rmk}\label{Onetoposmany_converso}
\emph{We note that the theorem implies that, given a site $({\cal C}, J)$, if all the objects of the form $l^{\cal C}_{J}(c)$ for $c\in {\cal C}$ are indecomposable objects (resp. atoms) of $\Sh({\cal C}, J)$ then the topos $\Sh({\cal C}, J)$ is locally connected (resp. atomic). Also, provided that $\cal C$ is cartesian, if all the objects of the form $l^{\cal C}_{J}(c)$ for $c\in {\cal C}$ are regular (resp. coherent) in $\Sh({\cal C}, J)$ then the topos $\Sh({\cal C}, J)$ is regular (resp. coherent). An application of this remark in the context of quotients of theories of presheaf type will be provided in section \ref{Onetoposmany_last} below.}   
\end{rmk}

\section{Syntactic criteria}

Let $\mathbb T$ be a cartesian (resp. regular, coherent, geometric) theory over a signature $\Sigma$. We denote by ${\cal C}_{\mathbb T}^{\textrm{cart}}$ (resp. ${\cal C}_{\mathbb T}^{\textrm{reg}}$, ${\cal C}_{\mathbb T}^{\textrm{coh}}$, ${\cal C}_{\mathbb T}^{\textrm{geom}}$) the cartesian (resp. regular, coherent, geometric) syntactic category of $\mathbb T$ and by $J^{\textrm{reg}}_{\mathbb T}$ the regular (resp. coherent, geometric) topology on ${\cal C}_{\mathbb T}^{\textrm{reg}}$ (resp. ${\cal C}_{\mathbb T}^{\textrm{coh}}$, ${\cal C}_{\mathbb T}^{\textrm{geom}}$). Recall from \cite{El2} that if $\mathbb T$ is cartesian (resp. regular, coherent, geometric) then $[{{\cal C}_{\mathbb T}^{\textrm{cart}}}^{\textrm{op}}, \Set]$ (resp. $\Sh({\cal C}_{\mathbb T}^{\textrm{reg}}, J^{\textrm{reg}}_{\mathbb T})$, $\Sh({\cal C}_{\mathbb T}^{\textrm{coh}}, J^{\textrm{coh}}_{\mathbb T})$, $\Sh({\cal C}_{\mathbb T}^{\textrm{geom}}, J^{\textrm{geom}}_{\mathbb T})$) is a classifying topos for $\mathbb T$. Let us denote by $y^{\textrm{cart}}_{\mathbb T}:{\cal C}^{\textrm{cart}}_{\mathbb T}\to [{{\cal C}_{\mathbb T}^{\textrm{cart}}}^{\textrm{op}}, \Set]$ (resp. $y^{\textrm{reg}}_{\mathbb T}:{\cal C}^{\textrm{reg}}_{\mathbb T}\to \Sh({\cal C}_{\mathbb T}^{\textrm{reg}}, J^{\textrm{reg}}_{\mathbb T})$, $y^{\textrm{coh}}_{\mathbb T}:{\cal C}^{\textrm{coh}}_{\mathbb T}\to \Sh({\cal C}_{\mathbb T}^{\textrm{coh}}, J^{\textrm{coh}}_{\mathbb T})$, $y^{\textrm{geom}}_{\mathbb T}:{\cal C}^{\textrm{geom}}_{\mathbb T}\to \Sh({\cal C}_{\mathbb T}^{\textrm{geom}}, J^{\textrm{geom}}_{\mathbb T})$) the Yoneda embeddings.  

Let us introduce the following notions. Below, by a $\mathbb T$-provably functional geometric formula from $\{\vec{x}. \phi\}$ to $\{\vec{y}. \psi\}$ we mean a geometric formula $\theta(\vec{x}, \vec{y})$ such that the sequents
$(\phi \vdash_{\vec{x}} (\exists \vec{y})\theta)$, $(\theta \vdash_{\vec{y}, \vec{x}} \phi \wedge \psi)$ and $(\theta \wedge \theta[\vec{z}\slash \vec{y}] \vdash_{\vec{x}, \vec{y}, \vec{z}} \vec{y}=\vec{z})$ are provable in $\mathbb T$.  

\begin{definition}
Let $\mathbb T$ be a geometric theory over a signature $\Sigma$ and $\phi(\vec{x})$ a geometric formula-in-context over $\Sigma$. Then

(i) we say that $\phi(\vec{x})$ is $\mathbb T$-complete if the sequent ($\phi \vdash_{\vec{x}} \bot$) is not provable in $\mathbb T$, and for every geometric formula $\phi$ in the same context either ($\chi \vdash_{\vec{x}} \bot$) or ($\chi \wedge \phi \vdash_{\vec{x}} \bot$) is provable in $\mathbb T$;

(ii) we say that $\phi(\vec{x})$ is $\mathbb T$-indecomposable if for any family $\{\psi_{i}(\vec{x}) \textrm{ | } i\in I\}$ of geometric formulae in the same context such that for each $i$, $\psi_{i}$ $\mathbb T$-provably implies $\phi$ and for any distinct $i, j\in I$, $\psi_{i}\wedge \psi_{j} \vdash_{\vec{x}} \bot$ is provable in $\mathbb T$, we have that $\phi \vdash_{\vec{x}} \mathbin{\mathop{\textrm{\huge $\vee$}}\limits_{i\in I}} \psi_{i}$ provable in $\mathbb T$ implies $\phi \vdash_{\vec{x}} \psi_{i}$ provable in $\mathbb T$ for some $i\in I$;

(iii) we say that $\phi(\vec{x})$ is $\mathbb T$-irreducible if for any family $\{\theta_{i} \textrm{ | } i\in I\}$ of $\mathbb T$-provably functional geometric formulae $\{\vec{x_{i}}, \vec{x}.\theta_{i}\}$ from $\{\vec{x_{i}}. \phi_{i}\}$ to $\{\vec{x}. \phi\}$ such that $\phi \vdash_{\vec{x}} \mathbin{\mathop{\textrm{\huge $\vee$}}\limits_{i\in I}}(\exists \vec{x_{i}})\theta_{i}$ is provable in $\mathbb T$, there exist $i\in I$ and a $\mathbb T$-provably functional geometric formula $\{\vec{x}, \vec{x_{i}}. \theta'\}$ from $\{\vec{x}. \phi\}$ to $\{\vec{x_{i}}. \phi_{i}\}$ such that $\phi \vdash_{\vec{x}} (\exists \vec{x_{i}})(\theta' \wedge \theta_{i})$ is provable in $\mathbb T$;

(iv) we say that $\phi(\vec{x})$ is $\mathbb T$-compact if for any family $\{\psi_{i}(\vec{x}) \textrm{ | } i\in I\}$ of geometric formulae in the same context, $\phi \vdash_{\vec{x}} \mathbin{\mathop{\textrm{\huge $\vee$}}\limits_{i\in I}} \psi_{i}$ provable in $\mathbb T$ implies $\phi \vdash_{\vec{x}} \mathbin{\mathop{\textrm{\huge $\vee$}}\limits_{i\in I'}} \psi_{i}$ provable in $\mathbb T$ for some finite subset $I'$ of $I$;

(v) we say that $\phi(\vec{x})$ is $\mathbb T$-supercompact if for any family $\{\psi_{i}(\vec{x}) \textrm{ | } i\in I\}$ of geometric formulae in the same context, $\phi \vdash_{\vec{x}} \mathbin{\mathop{\textrm{\huge $\vee$}}\limits_{i\in I}} \psi_{i}$ provable in $\mathbb T$ implies $\phi \vdash_{\vec{x}} \psi_{i}$ provable in $\mathbb T$ for some $i\in I$.   
\end{definition}

\begin{lemma}\label{Onetoposmany_equiv}
Let $\mathbb T$ be a geometric theory over a signature $\Sigma$ and $\phi(\vec{x})$ a geometric formula-in-context over $\Sigma$. Then

(i) $\phi(\vec{x})$ is $\mathbb T$-complete if and only if $y^{\textrm{geom}}_{\mathbb T}(\{\vec{x}. \phi\})$ is an atom of\\ $\Sh({\cal C}_{\mathbb T}^{\textrm{geom}}, J^{\textrm{geom}}_{\mathbb T})$;

(ii) $\phi(\vec{x})$ is $\mathbb T$-indecomposable if and only if $y^{\textrm{geom}}_{\mathbb T}(\{\vec{x}. \phi\})$ is an indecomposable object of $\Sh({\cal C}_{\mathbb T}^{\textrm{geom}}, J^{\textrm{geom}}_{\mathbb T})$;

(iii) $\phi(\vec{x})$ is $\mathbb T$-irreducible if and only if $y^{\textrm{geom}}_{\mathbb T}(\{\vec{x}. \phi\})$ is an irreducible object of $\Sh({\cal C}_{\mathbb T}^{\textrm{geom}}, J^{\textrm{geom}}_{\mathbb T})$;

(iv) $\phi(\vec{x})$ is $\mathbb T$-compact if and only if $y^{\textrm{geom}}_{\mathbb T}(\{\vec{x}. \phi\})$ is a compact object of $\Sh({\cal C}_{\mathbb T}^{\textrm{geom}}, J^{\textrm{geom}}_{\mathbb T})$;

(v) $\phi(\vec{x})$ is $\mathbb T$-supercompact if and only if $y^{\textrm{geom}}_{\mathbb T}(\{\vec{x}. \phi\})$ is a supercompact object of $\Sh({\cal C}_{\mathbb T}^{\textrm{geom}}, J^{\textrm{geom}}_{\mathbb T})$.     
\end{lemma}

\begin{proofs}
Parts (iii), (iv) and (v) follow immediately from Proposition \ref{Onetoposmany_express}, by using Lemma D1.4.4(iv) \cite{El2} and the fact that cover-mono factorizations of arrow exist in ${\cal C}_{\mathbb T}^{\textrm{geom}}$.

Parts (i) and (ii) follow from the fact that in a topos coproducts are the same thing as disjoint unions, since $y^{\textrm{geom}}_{\mathbb T}({\cal C}_{\mathbb T}^{\textrm{geom}})$ is closed in $\Sh({\cal C}_{\mathbb T}^{\textrm{geom}}, J^{\textrm{geom}}_{\mathbb T})$ under taking subobjects and $y^{\textrm{geom}}_{\mathbb T}:{\cal C}^{\textrm{geom}}_{\mathbb T}\to \Sh({\cal C}_{\mathbb T}^{\textrm{geom}}, J^{\textrm{geom}}_{\mathbb T})$ is a full and faithful geometric functor. 
\end{proofs}

Recall from \cite{El2} that a Grothendieck topos $\cal E$ is compact if and only if the terminal object of $\cal E$ is compact; hence, a geometric theory $\mathbb T$ over a signature $\Sigma$ is classified by a compact topos if and only if $y^{\textrm{geom}}_{\mathbb T}(\{[]. \top\})$ is compact in $\Sh({\cal C}_{\mathbb T}^{\textrm{geom}}, J^{\textrm{geom}}_{\mathbb T})$, if and only if for any family $\{\psi_{i} \textrm{ | } i\in I\}$ of geometric sentences, $\phi \vdash \mathbin{\mathop{\textrm{\huge $\vee$}}\limits_{i\in I}} \psi_{i}$ provable in $\mathbb T$ implies $\phi \vdash \mathbin{\mathop{\textrm{\huge $\vee$}}\limits_{i\in I'}} \psi_{i}$ provable in $\mathbb T$ for some finite subset $I'$ of $I$.

We note that any cartesian (resp. regular, coherent) theory $\mathbb T$ over a signature $\Sigma$ can be regarded as a geometric theory, and hence we have an equivalence of classifying toposes $[{{\cal C}_{\mathbb T}^{\textrm{cart}}}^{\textrm{op}}, \Set]\simeq \Sh({\cal C}_{\mathbb T}^{\textrm{geom}}, J^{\textrm{geom}}_{\mathbb T})$ (resp. $\Sh({\cal C}_{\mathbb T}^{\textrm{reg}}, J^{\textrm{reg}}_{\mathbb T}) \simeq \Sh({\cal C}_{\mathbb T}^{\textrm{geom}}, J^{\textrm{geom}}_{\mathbb T})$, $\Sh({\cal C}_{\mathbb T}^{\textrm{coh}}, J^{\textrm{coh}}_{\mathbb T}) \simeq \Sh({\cal C}_{\mathbb T}^{\textrm{geom}}, J^{\textrm{geom}}_{\mathbb T})$); moreover, it is immediate to see that for any cartesian (resp. regular, coherent) formula $\phi(\vec{x})$ over $\Sigma$, $y^{\textrm{cart}}_{\mathbb T}(\{\vec{x}. \phi\})$ (resp. $y^{\textrm{reg}}_{\mathbb T}(\{\vec{x}. \phi\})$, $y^{\textrm{coh}}_{\mathbb T}(\{\vec{x}. \phi\})$) corresponds under the equivalence to $y^{\textrm{geom}}_{\mathbb T}(\{\vec{x}. \phi\})$.

This remark, combined with Propositions \ref{Onetoposmany_identif} and \ref{Onetoposmany_express}, leads to the following results. Below, given a geometric theory $\mathbb T$ over a signature $\Sigma$, by saying that $\mathbb T$ is equivalent to a cartesian (resp. regular, coherent) theory we mean that $\mathbb T$ can be axiomatized by cartesian (resp. regular, coherent) sequents over $\Sigma$. 

\begin{theorem}\label{Onetoposmany_cartcrit}
Let $\mathbb T$ be a geometric theory over a signature $\Sigma$. Then $\mathbb T$ is equivalent to a cartesian theory if and only if for any cartesian formula $\{\vec{x}. \phi\}$ over $\Sigma$, for any family $\{\theta_{i} \textrm{ | } i\in I\}$ of $\mathbb T$-provably functional geometric formulae $\{\vec{x_{i}}, \vec{x}.\theta_{i}\}$ from $\{\vec{x_{i}}. \phi_{i}\}$ to $\{\vec{x}. \phi\}$ such that $\phi \vdash_{\vec{x}} \mathbin{\mathop{\textrm{\huge $\vee$}}\limits_{i\in I}}(\exists \vec{x_{i}})\theta_{i}$ is provable in $\mathbb T$, there exist $i\in I$ and a $\mathbb T$-provably functional geometric formula $\{\vec{x}, \vec{x_{i}}. \theta'\}$ from $\{\vec{x}. \phi\}$ to $\{\vec{x_{i}}. \phi_{i}\}$ such that $\phi \vdash_{\vec{x}} (\exists \vec{x_{i}})(\theta' \wedge \theta_{i})$ is provable in $\mathbb T$. 
\end{theorem}

\begin{proofs}
Let us suppose that $\mathbb T$ is cartesian. Then the property of $\{\vec{x}. \phi\}$ in the statement of the proposition is equivalent, by Lemma \ref{Onetoposmany_equiv}(iii), to saying that $y^{\textrm{geom}}_{\mathbb T}(\{\vec{x}. \phi\})$ is irreducible in $\Sh({\cal C}_{\mathbb T}^{\textrm{geom}}, J^{\textrm{geom}}_{\mathbb T})$. But this condition corresponds, under the equivalence $[{{\cal C}_{\mathbb T}^{\textrm{cart}}}^{\textrm{op}}, \Set]\simeq \Sh({\cal C}_{\mathbb T}^{\textrm{geom}}, J^{\textrm{geom}}_{\mathbb T})$, to saying that $y^{\textrm{cart}}_{\mathbb T}(\{\vec{x}. \phi\})$ is irreducible (equivalently, indecomposable projective) in $[{{\cal C}_{\mathbb T}^{\textrm{cart}}}^{\textrm{op}}, \Set]$, and this true (cfr. Proposition \ref{Onetoposmany_identif}(iii)).

Conversely, if $\mathbb T$ is geometric and for any cartesian formula $\phi(\vec{x})$ over $\Sigma$, $y^{\textrm{geom}}_{\mathbb T}(\{\vec{x}. \phi\})$ is irreducible in $\Sh({\cal C}_{\mathbb T}^{\textrm{geom}}, J^{\textrm{geom}}_{\mathbb T})$ then, denoted by $\cal G$ the full subcategory of ${\cal C}_{\mathbb T}^{\textrm{geom}}$ on the cartesian formulae, equivalently the cartesian syntactic category ${{\cal C}_{{\mathbb T}'}^{\textrm{cart}}}$ of the cartesianization ${\mathbb T}'$ of $\mathbb T$ (i.e. the theory axiomatized by all the cartesian sequents over $\Sigma$ which are provable in $\mathbb T$), we have that $\cal G$ is $J^{\textrm{geom}}_{\mathbb T}$-dense (by Lemma D1.3.8 \cite{El2}) and $J^{\textrm{geom}}_{\mathbb T}|_{\cal G}$ is the trivial Grothendieck topology. Thus the Comparison Lemma yields an equivalence $\Sh({\cal C}_{\mathbb T}^{\textrm{geom}}, J^{\textrm{geom}}_{\mathbb T})\simeq [{{\cal C}_{{\mathbb T}'}^{\textrm{cart}}}^{\textrm{op}}, \Set]$. Clearly, this equivalence sends the universal model of $\mathbb T$ in $\Sh({\cal C}_{\mathbb T}^{\textrm{geom}}, J^{\textrm{geom}}_{\mathbb T})$ to the universal model of ${\mathbb T}'$ in $[{{\cal C}_{{\mathbb T}'}^{\textrm{cart}}}^{\textrm{op}}, \Set]$, and hence, universal models being conservative, $\mathbb T$ and ${\mathbb T}'$ prove exactly the same geometric sequents over $\Sigma$ i.e. they are equivalent, as required.         
\end{proofs}

The `only if' direction in the following result extends Lemma D3.3.11 \cite{El2}.
\begin{theorem}\label{Onetoposmany_regcrit}
Let $\mathbb T$ be a geometric theory over a signature $\Sigma$. Then $\mathbb T$ is equivalent to a regular theory if and only if for any regular formula $\{\vec{x}. \phi\}$ over $\Sigma$, for any family $\{\psi_{i}(\vec{x}) \textrm{ | } i\in I\}$ of geometric formulae in the same context, $\phi \vdash_{\vec{x}} \mathbin{\mathop{\textrm{\huge $\vee$}}\limits_{i\in I}} \psi_{i}$ provable in $\mathbb T$ implies $\phi \vdash_{\vec{x}} \psi_{i}$ provable in $\mathbb T$ for some $i\in I$. 
\end{theorem}

\begin{proofs}
This follows similarly to Theorem \ref{Onetoposmany_cartcrit} by using Lemma \ref{Onetoposmany_equiv}(v) and Proposition \ref{Onetoposmany_identif}(iv).
\end{proofs}

\begin{theorem}\label{Onetoposmany_cohcrit}
Let $\mathbb T$ be a geometric theory over a signature $\Sigma$. Then $\mathbb T$ is equivalent to a coherent theory if and only if for any coherent formula $\{\vec{x}. \phi\}$ over $\Sigma$, for any family $\{\psi_{i}(\vec{x}) \textrm{ | } i\in I\}$ of geometric formulae in the same context, $\phi \vdash_{\vec{x}} \mathbin{\mathop{\textrm{\huge $\vee$}}\limits_{i\in I}} \psi_{i}$ provable in $\mathbb T$ implies $\phi \vdash_{\vec{x}} \mathbin{\mathop{\textrm{\huge $\vee$}}\limits_{i\in I'}} \psi_{i}$ provable in $\mathbb T$ for some finite subset $I'$ of $I$. 
\end{theorem}

\begin{proofs}
This follows similarly to Theorem \ref{Onetoposmany_cartcrit} by using Lemma \ref{Onetoposmany_equiv}(iv) and Proposition \ref{Onetoposmany_identif}(iv). 
\end{proofs}

\subsection{Locally connected theories}

The following result gives a syntactic characterization of the class of geometric theories classified by a locally connected (resp. connected and locally connected) topos.

\begin{theorem}\label{Onetoposmany_loccon}
Let $\mathbb T$ be a geometric theory over a signature $\Sigma$. Then $\mathbb T$ is classified by a locally connected topos (resp. connected and locally connected topos) if and only if for any geometric formula $\phi(\vec{x})$ over $\Sigma$ there exists a (unique) family $\{\psi_{i}(\vec{x}) \textrm{ | } i\in I\}$ of $\mathbb T$-indecomposable geometric formulae in the same context such that

(i) for each $i$, $\psi_{i}$ $\mathbb T$-provably implies $\phi$,

(ii) for any distinct $i, j\in I$, $\psi_{i}\wedge \psi_{j} \vdash_{\vec{x}} \bot$ is provable in $\mathbb T$, and

(iii) $\phi \vdash_{\vec{x}} \mathbin{\mathop{\textrm{\huge $\vee$}}\limits_{i\in I}} \psi_{i}$ is provable in $\mathbb T$\\ 
(resp. and $\{[].\top\}$ is $\mathbb T$-indecomposable). 
\end{theorem}

\begin{proofs}
Let us suppose that the classifying topos $\Sh({\cal C}_{\mathbb T}^{\textrm{geom}}, J^{\textrm{geom}}_{\mathbb T})$ of $\mathbb T$ is locally connected. Then, by Lemma D3.3.6 \cite{El2}, given a geometric formula $\phi(\vec{x})$ over $\Sigma$, $y^{\textrm{geom}}_{\mathbb T}(\{\vec{x}. \phi\})$ is uniquely expressible as a coproduct of indecomposable objects of $\Sh({\cal C}_{\mathbb T}^{\textrm{geom}}, J^{\textrm{geom}}_{\mathbb T})$. Now, since $y^{\textrm{geom}}_{\mathbb T}({\cal C}_{\mathbb T}^{\textrm{geom}})$ is closed in $\Sh({\cal C}_{\mathbb T}^{\textrm{geom}}, J^{\textrm{geom}}_{\mathbb T})$ under taking subobjects ($\Sh({\cal C}_{\mathbb T}^{\textrm{geom}}, J^{\textrm{geom}}_{\mathbb T})$ being the $\infty$-pretopos completion of ${\cal C}_{\mathbb T}^{\textrm{geom}}$), we can suppose that all the subobjects of $y^{\textrm{geom}}_{\mathbb T}(\{\vec{x}. \phi\})$ in $\Sh({\cal C}_{\mathbb T}^{\textrm{geom}}, J^{\textrm{geom}}_{\mathbb T})$, and in particular the indecomposable objects arising in our coproduct, are of the form $y^{\textrm{geom}}_{\mathbb T}(c)$ for some $c\in {\cal C}_{\mathbb T}^{\textrm{geom}}$. The condition of the criterion then follows from Lemma \ref{Onetoposmany_equiv}(ii) and the fact that the functor $y^{\textrm{geom}}_{\mathbb T}$ is geometric and full and faithful.

Conversely, if the condition of the criterion is satisfied then we have, by Lemma \ref{Onetoposmany_equiv}(ii) and Proposition \ref{Onetoposmany_denseness}(ii), that the objects of the form $y^{\textrm{geom}}_{\mathbb T}(\{\vec{y}. \psi\})$ for a $\mathbb T$-indecomposable formula $\psi(\vec{y})$ form a separating set for\\ $\Sh({\cal C}_{\mathbb T}^{\textrm{geom}}, J^{\textrm{geom}}_{\mathbb T})$ made of indecomposable objects; then $\Sh({\cal C}_{\mathbb T}^{\textrm{geom}}, J^{\textrm{geom}}_{\mathbb T})$ is locally connected by Theorem \ref{Onetoposmany_car}(i). 
\end{proofs}

The following result is the coherent analogue of this theorem.

\begin{theorem}
Let $\mathbb T$ be a coherent theory over a signature $\Sigma$. Then $\mathbb T$ is classified by a locally connected topos (resp. connected and locally connected topos) if and only if for any coherent formula $\phi(\vec{x})$ over $\Sigma$ there exists a (unique) finite family $\{\psi_{i}(\vec{x}) \textrm{ | } i\in I\}$ of $\mathbb T$-indecomposable geometric formulae in the same context such that

(i) for each $i$, $\psi_{i}$ $\mathbb T$-provably implies $\phi$,

(ii) for any distinct $i, j\in I$, $\psi_{i}\wedge \psi_{j} \vdash_{\vec{x}} \bot$ is provable in $\mathbb T$, and

(iii) $\phi \vdash_{\vec{x}} \mathbin{\mathop{\textrm{\huge $\vee$}}\limits_{i\in I}} \psi_{i}$ is provable in $\mathbb T$\\ 
(resp. and $\{[].\top\}$ is $\mathbb T$-indecomposable). 
\end{theorem}

\begin{proofs}
The proof proceeds analogously to the proof of Theorem \ref{Onetoposmany_loccon}, by using Theorem \ref{Onetoposmany_cohcrit}.
\end{proofs}

\begin{rmk}\label{Onetoposmany_weaker}
\emph{From the proof of the theorems it is clear that, by using the notion of dense subcategory, one can obtain alternative (although equivalent) versions of the criteria. For example, a weaker (in the `if' direction) version of the criterion of Theorem \ref{Onetoposmany_loccon} reads as follows: a geometric theory $\mathbb T$ is classified by a locally connected topos (resp. connected and locally connected topos) if and only if there exists a collection $\cal F$ (resp. a collection $\cal F$ containing $\{[].\top\}$) of $\mathbb T$-indecomposable geometric formulae-in-context over $\Sigma$ such that for any geometric formula $\{\vec{y}. \psi\}$ over $\Sigma$, there exist objects $\{\vec{x_{i}}. \phi_{i}\}$ in $\cal F$ as $i$ varies in $I$ and $\mathbb T$-provably functional geometric formulae $\{\vec{x_{i}}, \vec{y}.\theta_{i}\}$ from $\{\vec{x_{i}}. \phi_{i}\}$ to $\{\vec{y}. \psi\}$ such that $\psi \vdash_{\vec{y}} \mathbin{\mathop{\textrm{\huge $\vee$}}\limits_{i\in I}}(\exists \vec{x_{i}})\theta_{i}$ is provable in $\mathbb T$. Naturally, the `coherent' version of this criterion also holds.}  
\end{rmk}

We note that, by Theorem \ref{Onetoposmany_cohcrit}, if $\mathbb T$ is coherent and $\phi(\vec{x})$ is a coherent formula over $\Sigma$ then $\phi(\vec{x})$ is $\mathbb T$-indecomposable if and only if for any finite family $\{\psi_{i}(\vec{x}) \textrm{ | } i\in I\}$ of geometric formulae in the same context such that for each $i$, $\psi_{i}$ $\mathbb T$-provably implies $\phi$ and for any distinct $i, j\in I$, $\psi_{i}\wedge \psi_{j} \vdash_{\vec{x}} \bot$ is provable in $\mathbb T$, we have that $\phi \vdash_{\vec{x}} \mathbin{\mathop{\textrm{\huge $\vee$}}\limits_{i\in I}} \psi_{i}$ provable in $\mathbb T$ implies $\phi \vdash_{\vec{x}} \psi_{i}$ provable in $\mathbb T$ for some $i\in I$. 

Regarding regular theories, their classifying toposes are always connected and locally connected; this was already observed in \cite{El2}, and also follows from our Remark \ref{Onetoposmany_Rem}. 

We have already noticed that any atomic topos $\cal E$ is locally connected; in fact, as it is observed in \cite{El2}, the atoms of $\cal E$ are exactly the connected (equivalently, indecomposable) objects of $\cal E$. Hence, in view of Theorem \ref{Onetoposmany_car}, by replacing `$\mathbb T$-indecomposable' with `$\mathbb T$-complete' everywhere in the criteria above, one obtains syntactic criteria for a geometric (resp. coherent) theory to be classified by an atomic topos; also, one can obtain alternative versions of these criteria in the same spirit as in Remark \ref{Onetoposmany_weaker}.

Concerning atomic toposes, let us notice that if $\phi(\vec{x})$ is a $\mathbb T$-complete formula then $\phi(\vec{x})$ is $\mathbb T$-provably equivalent to a regular formula; indeed, this follows immediately from the fact that any geometric formula is provably equivalent to a disjunction of regular formulae (Lemma D1.3.8 \cite{El2}).

We can give the following criterion for a regular theory to be classified by an atomic topos.

\begin{proposition}
Let $\mathbb T$ be a regular theory over a signature $\Sigma$. Then $\mathbb T$ is classified by an atomic topos if and only if every regular formula over $\Sigma$ is either $\mathbb T$-provably equivalent to $\bot$ or $\mathbb T$-complete.
\end{proposition}

\begin{proofs}
Let $\phi(\vec{x})$ be a regular formula over $\Sigma$. If $\mathbb T$ is classified by an atomic topos then, by the discussion above, $\phi$ is expressible as a disjunction of $\mathbb T$-complete regular formulae; but this implies, by Theorem \ref{Onetoposmany_regcrit}, that either $\phi$ is $\mathbb T$-provably equivalent to $\bot$ or it is $\mathbb T$-provably equivalent to one of these formulae and hence $\mathbb T$-complete.

Conversely, if every regular formula over $\Sigma$ is either $\mathbb T$-provably equivalent to $\bot$ or $\mathbb T$-complete then, by Lemma \ref{Onetoposmany_equiv} and the fact that the set of objects of the form $y^{\textrm{reg}}_{\mathbb T}(\{\vec{x}. \phi\})$ for a regular formula $\phi(\vec{x})$ over $\Sigma$ form a separating set for the classifying topos $\Sh({\cal C}_{\mathbb T}^{\textrm{reg}}, J^{\textrm{reg}}_{\mathbb T})$ of $\mathbb T$, we have that $\Sh({\cal C}_{\mathbb T}^{\textrm{reg}}, J^{\textrm{reg}}_{\mathbb T})$ of $\mathbb T$ has a separating set of atoms and hence, by Theorem \ref{Onetoposmany_car}(iii), it is atomic, as required.    
\end{proofs}

In passing, we note an interesting property of theories classified by atomic toposes.

\begin{theorem}
Let $\mathbb T$ be a regular (resp. coherent) theory over a signature $\Sigma$ which is classified by an atomic (equivalently, Boolean) topos. Then every geometric formula over $\Sigma$ is $\mathbb T$-provably equivalent to a regular (resp. coherent) formula over $\Sigma$.
\end{theorem}

\begin{proofs}
Let $\phi(\vec{x})$ be a geometric formula over $\Sigma$. Then, the classifying topos of $\mathbb T$ being Boolean, $\top \vdash_{\vec{x}} \phi(\vec{x}) \vee \neg^{\mathbb T} \phi(\vec{x})$ is provable in $\mathbb T$ (where $\neg^{\mathbb T} \phi(\vec{x})$ denotes the pseudocomplementation of $\phi(\vec{x})$ in $\Sub_{{\cal C}_{\mathbb T}^{\textrm{geom}}}(\{\vec{x}. \top\})$ as in \cite{OC7}). But, by Theorem \ref{Onetoposmany_regcrit}, $\top(\vec{x})$ is $\mathbb T$-supercompact (resp. $\mathbb T$-compact) and, since by Lemma D1.3.8 \cite{El2} $\phi(\vec{x})$ is ($\mathbb T$-)provably equivalent to a disjunction of regular formulae, it thus follows that $\phi(\vec{x})$ is $\mathbb T$-provably equivalent to a single regular formula (resp. a finite disjunction of regular formulae), as required.
\end{proofs}

\subsection{Theories of presheaf type}

In this section we give a characterization of the class of geometric (resp. coherent, regular) theories classified by a presheaf topos.\\ We recall that a theory classified by a presheaf topos is said to be of presheaf type.

Below, for a subcanonical site $({\cal C}, J)$, we denote by $y:{\cal C}\to \Sh({\cal C}, J)$ the factorization through $\Sh({\cal C}, J)\hookrightarrow [{\cal C}^{\textrm{op}}, \Set]$ of the Yoneda embedding.

\begin{theorem}\label{Onetoposmany_presheaf}
Let $({\cal C}, J)$ be a subcanonical site such that $y({\cal C})$ is closed in $\Sh({\cal C}, J)$ under retracts. Then $\Sh({\cal C}, J)$ is equivalent to a presheaf topos if and only if $J$ is rigid.
\end{theorem} 
 
\begin{proofs}
The `if' direction follows at once from the Comparison Lemma. Let us then prove the `only if' direction. If ${\cal E}=\Sh({\cal C}, J)$ is equivalent to a presheaf topos then, by Lemma C2.2.20 \cite{El2}, $\cal E$ has a separating set of indecomposable projective objects. Now, suppose A is an indecomposable projective in $\cal E$. Then, as it is observed in the proof of Lemma C2.2.20 \cite{El2}, given any epimorphic family $\{f_{i}: B_{i} \to A \textrm{ | } i\in I\}$, at least one $f_{i}$ must be a split epimorphism; in particular $A$ is $J_{\cal E}$-irreducible, where $J_{\cal E}$ is the canonical coverage on $\cal E$. Hence, by taking as epimorphic family the collection of all the arrows in $\cal E$ from objects of the form $y(c)$ to $A$, we obtain that $A$ is a retract in $\cal E$ of an object of the form $y(c)$. Thus, by our hypotheses, $A$ is itself, up to isomorphism, of the form $y(c)$ for some $c\in {\cal C}$. Let us denote by ${\cal C}'$ the full subcategory of $\cal C$ on the objects $c$ such that $y(c)$ is indecomposable and projective in $\cal E$; then the objects in $y({\cal C}')$ form a separating set for $\cal E$. Thus, for any object $B$ of $\cal E$ the family of all the arrows in $\cal E$ from objects of the form $y(c)$ for $c\in {\cal C}'$ to $B$ generates a $J_{\cal E}$-covering sieve. But, $J$ being subcanonical, $J=J_{\cal E}|_{\cal C}$ (by Proposition C2.2.16 \cite{El2}) and hence for any object $c\in {\cal C}$ the collection of all arrows in $\cal C$ from objects of ${\cal C}'$ to $c$ is $J$-covering; so, by Proposition \ref{Onetoposmany_express}(iii), $J$ is rigid. 
\end{proofs}

\begin{rmk}\label{Onetoposmany_cauchy}
\emph{We note that, under the hypotheses of the theorem, if $\cal C$ is Cauchy-complete (in particular if $\cal C$ is cartesian) then $y({\cal C})$ is closed in ${\cal E} = \Sh({\cal C}, J)$ under retracts. Indeed, let $i:A\mono y(c)$, $r:y(c)\epi A$ be a retract of $A$ in $\cal E$ i.e. $r\circ i=1_{A}$. Then $i\circ r:y(c)\to y(c)$ is idempotent. Now, since $y$ is full and faithful, $i\circ r=y(e)$ for some idempotent $e:c\to c$ in $\cal C$. Since $\cal C$ is Cauchy complete then $e$ splits as $s\circ t$ where $t\circ s=1$. Then $y(s)$ and $y(t)$ form a retract of $A$ and hence, by the uniqueness up to isomorphism of the splitting of an idempotent in a category, it follows that $r$ is isomorphic to $y(t)$ and $i$ is isomorphic to $y(s)$, and in particular $A$ is isomorphic to $y(dom(s))$.}
\end{rmk}  
By Remark \ref{Onetoposmany_cauchy}, the regular (resp. coherent, geometric) syntactic sites for regular (resp. coherent, geometric) theories all satisfy the hypotheses of Theorem \ref{Onetoposmany_presheaf}. Thus we obtain the following results.

\begin{corollary}
Let $\mathbb T$ be a geometric theory over a signature $\Sigma$. Then $\mathbb T$ is of presheaf type if and only if there exists a collection $\cal F$ of geometric formulae-in-context over $\Sigma$ satisfying the following properties:

(1) for any geometric formula $\{\vec{y}. \psi\}$ over $\Sigma$, there exist objects $\{\vec{x_{i}}. \phi_{i}\}$ in $\cal F$ as $i$ varies in $I$ and $\mathbb T$-provably functional geometric formulae $\{\vec{x_{i}}, \vec{y}.\theta_{i}\}$ from $\{\vec{x_{i}}. \phi_{i}\}$ to $\{\vec{y}. \psi\}$ such that $\psi \vdash_{\vec{y}} \mathbin{\mathop{\textrm{\huge $\vee$}}\limits_{i\in I}}(\exists \vec{x_{i}})\theta_{i}$ is provable in $\mathbb T$;

(2) for any formula $\{\vec{x}. \phi\}$ in $\cal F$, for any family $\{\theta_{i} \textrm{ | } i\in I\}$ of $\mathbb T$-provably functional geometric formulae $\{\vec{x_{i}}, \vec{x}.\theta_{i}\}$ from $\{\vec{x_{i}}. \phi_{i}\}$ to $\{\vec{x}. \phi\}$ such that $\phi \vdash_{\vec{x}} \mathbin{\mathop{\textrm{\huge $\vee$}}\limits_{i\in I}}(\exists \vec{x_{i}})\theta_{i}$ is provable in $\mathbb T$, there exist $i\in I$ and a $\mathbb T$-provably functional geometric formula $\{\vec{x}, \vec{x_{i}}. \theta'\}$ from $\{\vec{x}. \phi\}$ to $\{\vec{x_{i}}. \phi_{i}\}$ such that $\phi \vdash_{\vec{x}} (\exists \vec{x_{i}})(\theta' \wedge \theta_{i})$ is provable in $\mathbb T$.  
\end{corollary}\qed

Note that condition (2) in the corollary says precisely that $\{\vec{x}. \phi\}$ is $\mathbb T$-irreducible; in particular, $\{\vec{x}. \phi\}$ is $\mathbb T$-supercompact i.e. for any family $\{ \{\vec{x}.\phi_{i}\} \textrm{ | } i\in I\}$ of geometric formulae in the same context which $\mathbb T$-provably imply $\{\vec{x}. \phi\}$ and such that $\phi \vdash_{\vec{x}} \mathbin{\mathop{\textrm{\huge $\vee$}}\limits_{i\in I}} \phi_{i}$ is provable in $\mathbb T$, there exists $i\in I$ such that $\phi_{i}$ and $\phi$ are $\mathbb T$-provably equivalent.

The following results are the coherent and regular analogues of this corollary.

\begin{corollary}
Let $\mathbb T$ be a coherent theory over a signature $\Sigma$. Then $\mathbb T$ is of presheaf type if and only if there exists a collection $\cal F$ of coherent formulae-in-context over $\Sigma$ satisfying the following properties:

(1) for any coherent formula $\{\vec{y}. \psi\}$ over $\Sigma$, there exists a finite number of objects $\{\vec{x_{i}}. \phi_{i}\}$ in $\cal F$ as $i$ varies in $I$ and $\mathbb T$-provably functional coherent formulae $\{\vec{x_{i}}, \vec{y}.\theta_{i}\}$ from $\{\vec{x_{i}}. \phi_{i}\}$ to $\{\vec{y}. \psi\}$ such that $\psi \vdash_{\vec{y}} \mathbin{\mathop{\textrm{\huge $\vee$}}\limits_{i\in I}}(\exists \vec{x_{i}})\theta_{i}$ is provable in $\mathbb T$;

(2) for any formula $\{\vec{x}. \phi\}$ in $\cal F$, for any finite family $\{\theta_{i} \textrm{ | } i\in I\}$ of $\mathbb T$-provably functional coherent formulae $\{\vec{x_{i}}, \vec{x}.\theta_{i}\}$ from $\{\vec{x_{i}}. \phi_{i}\}$ to $\{\vec{x}. \phi\}$ such that $\phi \vdash_{\vec{x}} \mathbin{\mathop{\textrm{\huge $\vee$}}\limits_{i\in I}}(\exists \vec{x_{i}})\theta_{i}$ is provable in $\mathbb T$, there exist $i\in I$ and a $\mathbb T$-provably functional coherent formula $\{\vec{x}, \vec{x_{i}}. \theta'\}$ from $\{\vec{x}. \phi\}$ to $\{\vec{x_{i}}. \phi_{i}\}$ such that $\phi \vdash_{\vec{x}} (\exists \vec{x_{i}})(\theta' \wedge \theta_{i})$ is provable in $\mathbb T$.  
\end{corollary}\qed

\begin{corollary}
Let $\mathbb T$ be a regular theory over a signature $\Sigma$. Then $\mathbb T$ is of presheaf type if and only if there exists a collection $\cal F$ of regular formulae-in-context over $\Sigma$ satisfying the following properties:

(1) for any regular formula $\{\vec{y}. \psi\}$ over $\Sigma$, there exists an object $\{\vec{x}. \phi\}$ in $\cal F$ and a $\mathbb T$-provably functional formula $\{\vec{x}, \vec{y}.\theta\}$ from $\{\vec{x}. \phi\}$ to $\{\vec{y}. \psi\}$ such that $\psi \vdash_{\vec{y}} (\exists \vec{x})\theta$ is provable in $\mathbb T$;

(2) for any formula $\{\vec{y}. \psi\}$ in $\cal F$, for any $\mathbb T$-provably functional regular formulae $\{\vec{x}, \vec{y}.\theta\}$ from $\{\vec{x}. \phi\}$ to $\{\vec{y}. \psi\}$ such that $\psi \vdash_{\vec{y}} (\exists \vec{x})\theta$ is provable in $\mathbb T$, there exist a $\mathbb T$-provably functional regular formula $\{\vec{y}, \vec{x}. \theta'\}$ from $\{\vec{y}. \psi\}$ to $\{\vec{x}. \phi\}$ such that $\psi \vdash_{\vec{y}} (\exists \vec{x})(\theta' \wedge \theta_{i})$ is provable in $\mathbb T$.  
\end{corollary}\qed
    
\section{Syntactic properties of quotients of theories of presheaf type}\label{Onetoposmany_last}

Let $\mathbb T$ be a theory of presheaf type over a signature $\Sigma$. Then, by choosing a canonical Morita-equivalence for $\mathbb T$ (as in \cite{OC7}), we have an equivalence of classifying toposes $[\textrm{f.p.} {\mathbb T}\textrm{-mod}(\Set), \Set]\simeq \Sh({\cal C}_{\mathbb T}^{\textrm{geom}}, J^{\textrm{geom}}_{\mathbb T})$ sending the universal model $M_{\mathbb T}$ of $\mathbb T$ in the topos $[\textrm{f.p.} {\mathbb T}\textrm{-mod}(\Set), \Set]$ (as in Theorem 3.1 \cite{OC7}) into the universal model $U_{\mathbb T}$ of $\mathbb T$ in $\Sh({\cal C}_{\mathbb T}^{\textrm{geom}}, J^{\textrm{geom}}_{\mathbb T})$ (as described in \cite{OC7}). In particular, if $M\in \textrm{f.p.} {\mathbb T}\textrm{-mod}(\Set)$ is a $\mathbb T$-model presented by a formula $\phi(\vec{x})$ over $\Sigma$ then, denoted by $y:\textrm{f.p.} {\mathbb T}\textrm{-mod}(\Set)^{\textrm{op}} \to [\textrm{f.p.} {\mathbb T}\textrm{-mod}(\Set), \Set]$ the Yoneda embedding, $y(M)$ is equal to $[[\vec{x}. \phi]]_{M_{\mathbb T}}$ and hence corresponds, under the equivalence\\ $[\textrm{f.p.} {\mathbb T}\textrm{-mod}(\Set), \Set]\simeq \Sh({\cal C}_{\mathbb T}^{\textrm{geom}}, J^{\textrm{geom}}_{\mathbb T})$, to the functor  
$y^{\textrm{geom}}_{\mathbb T}(\{\vec{x}. \phi\})=[[\vec{x}. \phi]]_{U_{\mathbb T}}$.

Now, if ${\mathbb T}'$ is a quotient of $\mathbb T$ then the subtopos of $\Sh({\cal C}_{\mathbb T}^{\textrm{geom}}, J^{\textrm{geom}}_{\mathbb T})$ corresponding to it via Theorem 3.6 \cite{OC6} transfers via the equivalence\\ $[\textrm{f.p.} {\mathbb T}\textrm{-mod}(\Set), \Set]\simeq \Sh({\cal C}_{\mathbb T}^{\textrm{geom}}, J^{\textrm{geom}}_{\mathbb T})$ to a subtopos\\ $\Sh(\textrm{f.p.} {\mathbb T}\textrm{-mod}(\Set)^{\textrm{op}}, J)\hookrightarrow [\textrm{f.p.} {\mathbb T}\textrm{-mod}(\Set), \Set]$ of $[\textrm{f.p.} {\mathbb T}\textrm{-mod}(\Set), \Set]$; the topology $J$ will be called the associated $\mathbb T$-topology of ${\mathbb T}'$. This gives rise to an equivalence $\Sh(\textrm{f.p.} {\mathbb T}\textrm{-mod}(\Set)^{\textrm{op}}, J) \simeq \Sh({\cal C}_{{\mathbb T}'}^{\textrm{geom}}, J^{\textrm{geom}}_{{\mathbb T}'})$ of classifying toposes for ${\mathbb T}'$ which sends $l^{\textrm{f.p.} {\mathbb T}\textrm{-mod}(\Set)^{\textrm{op}}}_{J}(M)$ to the functor $y^{\textrm{geom}}_{{\mathbb T}'}(\{\vec{x}. \phi\})=[[\vec{x}. \phi]]_{U_{{\mathbb T}'}}$.   

The following result provides a link between `geometrical' properties of $J$ and syntactic properties of ${\mathbb T}'$. 

\begin{theorem}\label{Onetoposmany_sem}
Let $\mathbb T$ be a theory of presheaf type over a signature $\Sigma$, ${\mathbb T}'$ be a quotient of $\mathbb T$ with associated $\mathbb T$-topology $J$ on $\textrm{f.p.} {\mathbb T}\textrm{-mod}(\Set)^{\textrm{op}}$ and $\phi(\vec{x})$ be a geometric formula over $\Sigma$ which presents a $\mathbb T$-model $M$. Then

(i) $\phi(\vec{x})$ is $\mathbb T$-irreducible; in particular, $\phi(\vec{x})$ is $\mathbb T$-provably equivalent to a regular formula;

(ii) if the site $(\textrm{f.p.} {\mathbb T}\textrm{-mod}(\Set)^{\textrm{op}}, J)$ is locally connected (for example when $\textrm{f.p.} {\mathbb T}\textrm{-mod}(\Set)^{\textrm{op}}$ satisfies the right Ore condition and every $J$-covering sieve is non-empty) then $\phi(\vec{x})$ is ${\mathbb T}'$-indecomposable;
  
(iii) if $(\textrm{f.p.} {\mathbb T}\textrm{-mod}(\Set)^{\textrm{op}}$ satisfies the right Ore condition and $J$ is the atomic topology on $(\textrm{f.p.} {\mathbb T}\textrm{-mod}(\Set)^{\textrm{op}}$ then $\phi(\vec{x})$ is ${\mathbb T}'$-complete;

(iv) if every $J$-covering sieve on $M$ contains a $J$-covering sieve generated by a finite family of morphisms (resp. by a single morphism) then $\phi(\vec{x})$ is ${\mathbb T}'$-compact (resp. ${\mathbb T}'$-supercompact).
\end{theorem}

\begin{proofs}
(i) By Lemma \ref{Onetoposmany_equiv}(iii), $\phi(\vec{x})$ is $\mathbb T$-irreducible if and only if $y^{\textrm{geom}}_{\mathbb T}(\{\vec{x}. \phi\})$ is an irreducible object of $\Sh({\cal C}_{\mathbb T}^{\textrm{geom}}, J^{\textrm{geom}}_{\mathbb T})$. But, by the discussion above, this is equivalent to saying that $y(M)$ is irreducible in $[\textrm{f.p.} {\mathbb T}\textrm{-mod}(\Set), \Set]$, and this is is true by Proposition \ref{Onetoposmany_identif}(iii). The fact that $\phi(\vec{x})$ is $\mathbb T$-provably equivalent to a regular formula then follows from Lemma D1.3.8(ii) \cite{El2}.

(ii) and (iii) By Lemma \ref{Onetoposmany_equiv}(ii) (resp. Lemma \ref{Onetoposmany_equiv}(i)), $\phi(\vec{x})$ is ${\mathbb T}'$-indecom-\\posable (resp. ${\mathbb T}'$-complete) if and only if $y^{\textrm{geom}}_{{\mathbb T}'}(\{\vec{x}. \phi\})$ is an indecomposable object (resp. an atom) of $\Sh({\cal C}_{{\mathbb T}'}^{\textrm{geom}}, J^{\textrm{geom}}_{{\mathbb T}'})$; but this is equivalent to saying that $l^{{\textrm{f.p.} {\mathbb T}\textrm{-mod}(\Set)}^{\textrm{op}}}_{J}(M)$ is an indecomposable object (resp. an atom) of $\Sh(\textrm{f.p.} {\mathbb T}\textrm{-mod}(\Set)^{\textrm{op}}, J)$, and this is true by Proposition \ref{Onetoposmany_identif}(i) (resp. Proposition \ref{Onetoposmany_identif}(ii)).

(iv) By Lemma \ref{Onetoposmany_equiv}(iii), $\phi(\vec{x})$ is ${\mathbb T}'$-compact (resp. ${\mathbb T}'$-supercompact) if and only if $y^{\textrm{geom}}_{{\mathbb T}'}(\{\vec{x}. \phi\})$ is a compact (resp. supercompact) object of $\Sh({\cal C}_{{\mathbb T}'}^{\textrm{geom}}, J^{\textrm{geom}}_{{\mathbb T}'})$; but this is equivalent to saying that $l^{{\textrm{f.p.} {\mathbb T}\textrm{-mod}(\Set)}^{\textrm{op}}}_{J}(M)$ is a compact (resp. supercompact) object of $\Sh(\textrm{f.p.} {\mathbb T}\textrm{-mod}(\Set)^{\textrm{op}}, J)$, and this is true by Proposition \ref{Onetoposmany_express}(i) (resp. Proposition \ref{Onetoposmany_express}(ii)).    
\end{proofs}

\begin{rmk}
\emph{The theorem can be profitably applied in the context of cartesian theories. Indeed, if $\mathbb T$ is a cartesian theory then every cartesian formula over $\Sigma$ presents a $\mathbb T$-model so that we have an equivalence between $\textrm{f.p.} {\mathbb T}\textrm{-mod}(\Set)^{\textrm{op}}$ and the cartesian syntactic category of $\mathbb T$. Thus the theorem provides syntactic properties of cartesian formulae in particular quotients of $\mathbb T$; for example, for any such formula $\phi(\vec{x})$, part (iii) of the lemma yields that $\phi(\vec{x})$ is ${\mathbb T}'$-complete where ${\mathbb T}'$ is the Booleanization of $\mathbb T$ (as defined in \cite{OC3}).}
\end{rmk}

As an application of the notion of irreducible object in a topos, we can prove the following result.

\begin{theorem}\label{Onetoposmany_presheafcomplete}
Let $\mathbb T$ be a theory of presheaf type over a signature $\Sigma$. Then

(i) Any finitely presentable $\mathbb T$-model in $\Set$ is presented by a $\mathbb T$-irreducible geometric formula $\phi(\vec{x})$ over $\Sigma$;

(ii) Conversely, any $\mathbb T$-irreducible geometric formula $\phi(\vec{x})$ over $\Sigma$ presents a $\mathbb T$-model.

In particular, the category $\textrm{f.p.} {\mathbb T}\textrm{-mod}(\Set)^{\textrm{op}}$ is equivalent to the full subcategory of ${\cal C}_{\mathbb T}^{\textrm{geom}}$ on the $\mathbb T$-irreducible formulae. 
\end{theorem}

\begin{proofs}
We have already observed that we have an equivalence\\ $\tau:\Sh({\cal C}_{\mathbb T}^{\textrm{geom}}, J^{\textrm{geom}}_{\mathbb T}) \simeq [\textrm{f.p.} {\mathbb T}\textrm{-mod}(\Set), \Set]$ of classifying toposes for $\mathbb T$. Now, if ${\cal C}_{\mathbb T}^{\textrm{irr}}$ is the full subcategory of ${\cal C}_{\mathbb T}^{\textrm{geom}}$ on the $\mathbb T$-irreducible formulae then, by Theorem \ref{Onetoposmany_presheaf}, we have that $\Sh({\cal C}_{\mathbb T}^{\textrm{geom}}, J^{\textrm{geom}}_{\mathbb T}) \simeq [({\cal C}_{\mathbb T}^{\textrm{irr}})^{\textrm{op}}, \Set]$ via the Comparison Lemma (cfr. the proof of Theorem \ref{Onetoposmany_identif}). Now, if $\tilde{{\cal C}_{\mathbb T}^{\textrm{irr}}}\hookrightarrow {\cal C}_{\mathbb T}^{\textrm{geom}}$ is the Cauchy-completion of ${\cal C}_{\mathbb T}^{\textrm{irr}}$ then $[({\cal C}_{\mathbb T}^{\textrm{irr}})^{\textrm{op}}, \Set] \simeq [(\tilde{{\cal C}_{\mathbb T}^{\textrm{irr}}})^{\textrm{op}}, \Set]$ and the resulting equivalence $[\tilde{{\cal C}_{\mathbb T}^{\textrm{irr}}}^{\textrm{op}}, \Set]\simeq [\textrm{f.p.} {\mathbb T}\textrm{-mod}(\Set), \Set]$ restricts to an equivalence $l:\tilde{{\cal C}_{\mathbb T}^{\textrm{irr}}} \simeq \textrm{f.p.} {\mathbb T}\textrm{-mod}(\Set)$ between the subcategories of indecomposable projective objects. Now, given $\{\vec{x}.\phi\}\in \tilde{{\cal C}_{\mathbb T}^{\textrm{irr}}}$, $\tau$ sends the functor $y^{\textrm{geom}}_{\mathbb T}(\{\vec{x}. \phi\})=[[\vec{x}. \phi]]_{U_{\mathbb T}}$ to $y(l(\{\vec{x}. \phi\}))=[[\vec{x}. \phi]]_{M_{\mathbb T}}$, from which it follows that the model $l(\{\vec{x}. \phi\})$ is finitely presented by $\{\vec{x}. \phi\}$. Then, by Theorem \ref{Onetoposmany_sem}(i), $\phi(\vec{x})$ is $\mathbb T$-irreducible. So we conclude that $\tilde{{\cal C}_{\mathbb T}^{\textrm{irr}}}$ is equal to ${\cal C}_{\mathbb T}^{\textrm{irr}}$ i.e. ${\cal C}_{\mathbb T}^{\textrm{irr}}$ is Cauchy complete, and hence $l$ gives an equivalence $\tilde{{\cal C}_{\mathbb T}^{\textrm{irr}}} \simeq \textrm{f.p.} {\mathbb T}\textrm{-mod}(\Set)$. It is easy to verify that this equivalence coincide with the dualizing functor $d$ of Theorem 3.6 \cite{OC7}.      
\end{proofs}
    
\begin{rmk}    
\emph{As an application of Theorem \ref{Onetoposmany_presheafcomplete} and Remark \ref{Onetoposmany_converso}, suppose that $J$ is the associated $\mathbb T$-topology of a quotient ${\mathbb T}'$ of $\mathbb T$. Then, if for any $\mathbb T$-irreducible formula (equivalently, formula presenting a $\mathbb T$-model) $\phi(\vec{x})$, $\phi(\vec{x})$ is ${\mathbb T}'$-indecomposable (resp. ${\mathbb T}'$-complete) then the classifying topos $\Sh(\textrm{f.p.} {\mathbb T}\textrm{-mod}(\Set)^{\textrm{op}}, J)$ of ${\mathbb T}'$ is locally connected (resp. atomic); indeed, as observed above, $\phi(\vec{x})$ is ${\mathbb T}'$-indecomposable (resp. ${\mathbb T}'$-atomic) if and only if $y^{\textrm{geom}}_{{\mathbb T}'}(\{\vec{x}. \phi\})$ is an indecomposable object (resp. an atom) of the topos $\Sh({\cal C}_{{\mathbb T}'}^{\textrm{geom}}, J^{\textrm{geom}}_{{\mathbb T}'})$, if and only if $l^{{\textrm{f.p.} {\mathbb T}\textrm{-mod}(\Set)}^{\textrm{op}}}_{J}(M)$ is an indecomposable object (resp. an atom) of $\Sh(\textrm{f.p.} {\mathbb T}\textrm{-mod}(\Set)^{\textrm{op}}, J)$.} 
\end{rmk}    

\vspace{10 mm}
\begin{flushleft}
{\bf Acknowledgements:} I am grateful to my Ph.D. supervisor Peter Johnstone for useful discussions. 
\end{flushleft}

\end{document}